\documentclass[twoside]{article}
\begin{document}
\mathchardef\sym="3218
\markboth{IMRE MAJOR}{G-SMALE THEOREM}
\author{Imre Major}
\title{A G-version of Smale's theorem}
\date{}
\maketitle

\noindent A B S T R A C T $\vspace{.3cm}$

We will prove the equivariant version of Smale's transversality theorem: suppose that the compact Lie-group $G$ acts on the compact differentiable
manifold $M$ on which an invariant Morse-function $f$ and an invariant vector field $X$ are given so that $X$ is gradient-like with respect
to $f$ (i.e. $X(f)<0$ away from critical orbits and $X$ is the gradient of $f$ (w.r.t. a fixed invariant Riemannian metric) on some invariant
open subsets about critical orbits of $f$.) Given a bound $\varepsilon>0$ we will prove the existence of an invariant vector field $Y$ of class
$C^1$ for which vector field $X+Y$ is also gradient-like such that:\vspace{.4cm}

(a) $\|Y\|_1<\varepsilon\hspace{.3cm}(\|.\|_1$ is the $C^1$ norm).\vspace{.2cm}

(b) The intersection of the stable and unstable sets of vector field $X+Y$ taken

at a pair of critical orbits of $f$ is transverse when restricted to an orbit type of

the action.\vspace{1cm}

{\noindent\bf Key words:} Group actions, transversality.\vspace{.6cm}

\section{Introduction, basic concepts}

Suppose, that compact Lie group $G$ acts on the compact orientable smooth manifold $M^m$ and let $f:M\to{\bf R}$ be an invariant function (i.e.
$f(gx)=f(x)\hspace{.3cm}(g\in G)$). Fix an invariant metric $\bf g$.
If an orbit contains a critical point of $f$ then, by invariance, all points of this orbit are critical points
and the orbit itself is called a {\it critical orbit} of $f$. Fix an invariant Riemannian metric $\bf g$ and at a point $p\in M$ let$$\perp_p
:=[T_pG(p)]^\perp$$be the perpendicular complement of the tangent space to the orbit through point $p$ and let$$U_p:={\rm exp}_p(\perp_p
(\epsilon))$$where $\perp_p(\epsilon)$ is the $\epsilon$-disk about the origin in $\perp_p$ on which the exponential map of the Levi-Civita
connection of metric $\bf g$ is injective. Then a critical point $x$ of $f$ is also a critical point of $f|_{U_x}$. We say that the $G
$-invariant function $f$ is a {\it $G$-Morse function} if the Hessian of $f|_{U_x}$ at $x$ is non-degenerate for each critical point $x$.
This property does not depend on the choice of metric $\bf g$
(see e.g. Wasserman \cite{W}). Non-degeneracy of the Hessian also ensures that each critical orbit has an invariant neighborhood (called {\it tube}
about the orbit) that does not contain any other critical orbits. We can suppose, that $GU_x$ is such an invariant neighborhood.\vspace{.4cm}

\noindent The induced action of the {\it isotropy subgroup} $G_p$ at a point $p\in M$:$$G_p\times T_pM/T_pG(p)\to T_pM/T_pG(p)$$on the normal space is
called the {\it normal action} (see e.g. Bredon \cite{B}).$$pr:M\to M/G,\hspace{.7cm}p\longrightarrow G(p)$$ is the {\it canonical projection}.
For a set $N\subset M$ we use notation ${\cal N}:=pr(N)$.\vspace{.4cm}

{\noindent}Observe that the relation $$x\sim y\Leftrightarrow\exists(g\in G)\ \ni\ G_x=gG_yg^{-1}\ \ \ (x,y\in M)$$ is an equivalence relation.
The equivalence classes provide a partition\begin{eqnarray}M=\mathop{\cup}\limits_{\alpha\in{\cal A}} M_\alpha\end{eqnarray} of $M$. The index
set ${\cal A}$ is the set of conjugacy classes of isotropy subgroups of $G$. It is partially ordered by relation$$\alpha\prec\beta\ \Leftrightarrow\
\forall(x\in M_\alpha)\ \exists\ y\in M_\beta\ \ni\ G_y\subset G_x$$(the property on the RHS does not depend on the choice of representative
$x\in M_\alpha$).

An index $\alpha\in{\cal A}$ (or submanifold $M_\alpha$) is called an {\it orbit-type} and we say that a point $p$  is {\it of type} $\alpha$
if $p\in M_\alpha$ (in notation, $[p]=\alpha$). Note, that $M_\alpha\subset M$ is a $G$-invariant subset. Theorem 4. of Chapter 3. in \cite{B}
implies that $M_\alpha$ and ${\cal M}_\alpha$ are smooth manifolds and the partition in formula (1) is locally finite, thus $\cal A$ is finite
when $M$ is compact (cf. 4.3 Theorem, pg 187 and 10.4 Theorem, pg. 220 in \cite{B}). Let $m_\alpha:=$dim(${\cal M}_\alpha),\ m_\alpha^\perp:=$
codim$(M_\alpha)$ and $o_\alpha:=$dim$(G/G_p)\ (p\in M_\alpha$). Notice that then dim$(M):=m=m_\alpha+m_\alpha^\perp+o_\alpha$. An $\alpha
$-{\it slice} is a disk $D\subset M_\alpha$ which intersects an orbit at most once and along which the isotropy subgroup is constant, moreover
the union of orbits $GD$ is open in $M_\alpha$.

For an arbitrary subset $Q\subset M$ $$Q_\alpha:=Q\cap M_\alpha$$is the $\alpha${\it-part of} $Q$.\vspace{.4cm}

The set of differentials of left translations $\{{\rm d}L_g\ |\ g\in G\}$ act on the tangent bundle $TM$. A vector field $X$ is invariant under
this action (called an {\it invariant vector field}) {\bf iff} its flow is an eqivariant flow (i.e. for trajectory $\lambda_p$ of vector field
$X$ through point $p$ $$L_g(\lambda_p)=\lambda_{gp}$$holds.) This implies that the isotropy group is constant along trajectories, in particular,
an invariant vector field is tangent to the orbit types.\vspace{.4cm}

{\noindent\bf Definition 1.} An invariant vector field $X$ is {\it gradient-like} for the $G$-Morse function $f$ if:\vspace{.2cm}

(i) Each critical orbit has an invariant neighborhood $U$ such that$$X|_U=-{\rm grad}_{\bf g}(f)|_U$$

(ii) $X(f)<0$ away from critical orbits.\pagebreak

{\noindent\bf Definition 2.} At a critical orbit $O$:\vspace{.2cm}

$$W^-_O=\{p\in M|\mathop{\lim}\limits_{t\to -\infty}\lambda_p(t)\in O\}$$is called the {\it  unstable set} and

$$W^+_O=\{p\in M|\mathop{\lim}\limits_{t\to +\infty}\lambda_p(t)\in O\}$$is the {\it stable set} of the flow (of invariant vector field
$X$.) (Notations

\noindent $W^s_O=W^+_O,\ W_O^u=W_O^-$ are also used in the literature).\vspace{.4cm}

{\noindent\bf Definition 3.} A {\it Morse chart} about critical orbit $O$ is given by:

(i) A splitting of the normal bundle ${\bf\perp}_O$ of $O$ into two invariant orthogonal

subbundles ${\bf \perp}_O={\bf \perp}_O^-\oplus{\bf \perp}_O^+$.

(ii) An equivariant diffeomorphism $\eta_O:{\bf\perp}_O(\epsilon)\to U_O$ from the $\epsilon$-disc bundle of

the normal bundle of $O$ (with constant $\epsilon$) onto an invariant open neighborhood

$U_O$ of $O$ ($\eta_O$ is the identity on the zero section $O$ of ${\bf\perp}_O$) such that$$f\circ\eta_O=-\|P^-\|^2+\|P^+\|^2+f(O)$$

where $(P^-,P^+):{\bf\perp}_O\to ({\bf\perp}_O^-,{\bf\perp}_O^+)$ are the projections that belong to the

decomposition in (i). The open set $U_O$ is called a {\it Morse-tube}.\vspace{.2cm}

{\noindent\it Observe} that ${\bf\perp}_O^-\longrightarrow O$ and ${\bf\perp}_O^+\longrightarrow O$ are $G$-vector bundles induced from the
orthogonal representations on the Eucledian spaces $\perp_x^-:=\perp_x\cap{\bf\perp}_O^-,\ \perp_x^+:=\perp_x\cap{\bf\perp}_O^+$. The restriction
${\bf g}|_{\perp_O}$ is a scalar product on vector bundle $\perp_O$, thus it defines a Reimannian metric $\langle,\rangle$ on $\perp_O$ in the
canonical way. The push-forward $\eta_{O*}\langle,\rangle$ of this Riemannian metric along map $\eta_O$ can be patched together with original
metric $\bf g$ by using an invariant cutoff function. Thus we can presume that for the restrictions we have:$${\bf g}|_{U_O}=\eta_{O*}\langle,\rangle|_{U_O}$$

{\noindent\bf Lemma:} {\it(Equivariant Morse Lemma)} Let $f:M\to\bf R$ be a $G$-Morse function on a Riemannian $G$-manifold. Then there is a
Morse chart about each critical orbit (see Wasserman \cite{W}).\vspace{.2cm}

{\noindent\it Note} that by the above lemma the stable and unstable sets are, in fact, invariant submanifolds of $M$. Let $O_1,\dots,O_K$
be the set of critical orbits of $f$, and abbreviate $W_j^+:=W_{O_j}^+$, etc. \vspace{.2cm}

{\noindent\bf Definition 4.} The gradient-like vector field $X$ (or its flow $\Lambda:M\times{\bf R}\to M$) is $G${\it-Morse-Smale} if it is of class $C^1$
moreover $W^+_j\cap M_\alpha$ and $W^-_k\cap M_\alpha$ intersect transversely as submanifolds of $M_\alpha$ for each choice of $\alpha\in{
\cal A},\ 1\leq j,k\leq K$. (We refer to this property as {\it relative transversality} of stable and unstable submanifolds, or $\alpha$-transversality,
when the orbit type $\alpha$ is fixed.)\vspace{.2cm}

As we wish to perturb a given gradient-like vector field by an invariant vector field (the flow of which thus keeps orbit types), the above
definition seems to be the only plausible one for the $G$-version of Morse-Smale property.\vspace{.4cm}

{\noindent\bf Theorem:} {\it For any given $\varepsilon>0$ an invariant gradient-like vector field $X$ can be approximated by a G-Morse-Smale vector
field $X'$, which is also gradient-like for $f$ such that} $\|X-X'\|_1<\varepsilon.$\vspace{.5cm}

\section{Proof of the Theorem}

Fix an invariant gradient-like vector field $X$ for $G$-Morse function $f$ on Riemannian $G$-manifold $(M,{\bf g})$ (e.g. -grad$_{\bf g}(f)$
is such a vector field). In order to ensure that the perturbed vector field is $C^1$-close to the original one, we need a family of special
coordinate charts. Invariant charts would serve our purpose the best, however, we might not be able to arrange such charts with invariant
domains (by compactness of group $G$, the domain of such a chart can't be an $m$-ball).\vspace{.3cm}

{\noindent\it Remark:} As we will perturb vector field $X$ within the Morse-tubes (which we have already fixed), we could measure the $C^1$
norm of the perturbing vector field with respect to the Morse-charts (and in some sense, this is what will happen). Still, for completeness of our
discussion we wish to include the passage below, in which we introduce a special {\bf finite} set of charts that cover $M$.

As the boundary of an orbit type $M_\alpha$ is the union of some lower dimensional orbit types that preceed $\alpha$ with respect to $\prec$,
a {\it level} can be associated to each orbit type: the closed ones being at level $0$ and inductively, at level $i>0$ we have the orbit types such
that their boundaries contain orbit types of level at most $i-1$ and contains at least one orbit type at level $i-1$. (The level is also called
the {\it depth} in the literature.)

For each orbit type fix a representative isotropy subgroup $G_\alpha=G_{p_\alpha}$ (where $p_\alpha\in M_\alpha$), a coordinate chart $(E
_\alpha,(\tilde x_\alpha^{m_\alpha+m_\alpha^\perp+1},\dots,\tilde x_\alpha^m))$ about the unit element $G_\alpha$ of quotient group $G/G
_\alpha$, an open neighborhood $E_\alpha'\subset G/G_\alpha$ of $G_\alpha$ such that $\overline{E_\alpha'}\subset E_\alpha$ and elements
$g_\alpha^1,\dots,g_\alpha^{n_\alpha}\in G$ such that the translates$$g_\alpha^1E'_\alpha,\dots,g_\alpha^{n_\alpha}E'_\alpha$$cover $G/G
_\alpha$. Note that each orbit of type $\alpha$ contains a point with isotropy subgroup $G_\alpha$. Let ${\bf N}_\alpha\mathop{\longrightarrow}
\limits^{\Pi_\alpha}M_\alpha$ be the normal bundle of orbit type $M_\alpha$.

It is a well known fact that $G$ is a principal bundle over quotient group $G/G_\alpha$ by the canonical projection $G\longrightarrow
G/G_\alpha$ with structure group $G_\alpha$, so we can fix a section $\sigma_\alpha:E_\alpha\to G,\ \sigma_\alpha(G_\alpha)=e$ of this bundle over
contractible neighborhood $E_\alpha$ (i.e. $\sigma_\alpha(gG_\alpha)G_\alpha=gG_\alpha$).\vspace{.3cm}

{\noindent\bf Definition 5.} A coordinate chart $(U,(x^1,\dots,x^m))$ is {\it adapted to orbit type $\alpha$} if there exists a relatively
compact $\alpha$-slice $U^*_\alpha\subset M_\alpha$ with isotropy subgroup $G_\alpha$, a coordinate chart $({\cal U}^*_\alpha,x_\alpha^1,
\dots,x_\alpha^{m_\alpha}))$ (with respect to the smoothness structure on quotient manifold ${\cal M}_\alpha$), an index $1\leq i\leq
n_\alpha$ and $\epsilon>0$ such that:\vspace{.3cm}

(i) $U={\rm exp}({\bf N}_\alpha(\epsilon)|_{U_\alpha^i})$ with $U_\alpha^i=g_\alpha^iE_\alpha U^*_\alpha$, i.e. $U$ is the exp-image of the

restriction of the $\epsilon$-disc bundle of ${\bf N}_\alpha$ to subset $U^i_\alpha$.\vspace{.2cm}

(ii) $x^i(q)=x^i_\alpha\circ pr\circ\Pi_\alpha\circ{\rm exp}^{-1}(q)\hspace{.4cm}(i=1,\dots,m_\alpha,\ q\in U)$.\vspace{.2cm}

(iii) By clause (i), $\forall\ r\in U\ \exists!$ pair $(gG_\alpha,q)\in E_\alpha\times U^*_\alpha$ such that:$$\Pi_\alpha\circ{\rm exp}
^{-1}(r)=g_\alpha^igq$$

holds. Define \begin{eqnarray}x^i(r)=\tilde x_\alpha^i(gG_\alpha),\hspace{.8cm}(i=m_\alpha+m_\alpha^\perp+1,\dots,m)\end{eqnarray}\vspace{.2cm}

(iv) Fix an orthonormal frame bundle $({\bf v}_1,\dots,{\bf v}_{m^\perp_\alpha})$ of trivial bundle ${\bf N}_\alpha|_{U^*_\alpha}$

and extend it to $U_\alpha^i$ by ${\bf v}_j(g_\alpha^igq)={\rm d}L_{g_\alpha^i\sigma_\alpha(gG_\alpha)}({\bf v}_j(q))$. Define:\begin{eqnarray}{\exp}^{-1}(r)
=\mathop{\sum}\limits_{i=1}^{m^\perp_\alpha}x^{m_\alpha+i}(r){\bf v}_i\hspace{.8cm}(r\in U)\end{eqnarray}\vspace{.2cm}

\noindent Note that the first $m_\alpha$ coordinates do not depend on the choice of the rest of the coordinates. For the definition of the
$C^1$-norm of a vector field we need to fix a compact subset within each chart, which we can get as follows:\vspace{.2cm}

Fix $K^*_\alpha\subset U^*_\alpha$ compact set, $0<\epsilon'<\epsilon$ real number and define$$K={\rm exp}(\overline{\bf N}_\alpha(\epsilon')|
_{K_\alpha^i})\subset U$$

with $K_\alpha^i=g_\alpha^i\overline{E'_\alpha}K^*_\alpha$.\vspace{.4cm}

\noindent We will call int($K)$ the {\it strong interior} of chart $(U,(x^1,\dots,x^m))$. In the sequel we will presume that each adapted
chart has a fixed compact set in its interior (even if we don't state this explicitely).\vspace{.4cm}

\noindent To cover $M$ choose adapted charts$$(U_j,(x_j^1,\dots,x_j^m))\hspace{.7cm}(j=1,\dots,j_0)$$ so that their strong interior cover the
level-$0$ orbit types. Then the complement of the union of strong interiors of these charts in any of the level-$1$ strata is compact, so it can be
covered by the strong interiors of$$(U_j,(x_j^1,\dots,x_j^m))\hspace{.7cm}(j=j_0+1,\dots,j_1)$$so that each chart is adapted to the
stratum at issue, a.s.o. Finally we get a finite family of adapted coordinate charts so that their strong interiors int$(K_1),\dots,$int($K
_{j_L}$) cover $M$ (here $L$ is the highest level).\vspace{.3cm}

\noindent Given a vector field $Y$ with local coordinates:$$Y|_{U_j}=\mathop{\sum}\limits_{i=1}^my_j^i{\partial\over\partial x_j^i}\hspace
{.7cm}(j=1,\dots,j_L)$$its $C^1$-norm is defined as:\begin{eqnarray}\|Y\|_1:=\mathop{\sum}\limits_{i,j,k}\left(\mathop{\sup}\limits_{K_j}
|y_j^i|+\mathop{\sup}\limits_{K_j}|{\partial\over\partial x_j^k}y_j^i|\right)\end{eqnarray}\vspace{.1cm}

We will perturb vector field $X$ into a $G$-Morse-Smale vector field in succession of increasing order of the level of strata. Although
domains $U_1,\dots,U_{j_L}$ cover $M$, each stratum has parts covered by some of the $U_j$'s that are adapted to a stratum at a lower level.
To attain $\alpha$-transversality, we should adjust $X$ in charts that are adapted to orbit type $\alpha$. It seems so that in
order to end up with a $C^1$ vector field, we need to modify vector field $X$ in finite steps. To cover a non-compact
stratum, however, we need to use infinitely many adapted charts. To overcome this discrepancy, for each orbit type $\alpha$ we will choose a
countable family of adapted coordinate charts with set of domains $\underline Q^{(\alpha)}$ which is a finite union $$\underline Q^{(\alpha)}
=\underline Q^{(\alpha);1}\cup\dots\cup\underline Q^{(\alpha);k_\alpha}$$of sub-families such that:\vspace{.3cm}

{\bf Property I.} For a fixed $i=1,\dots,k_\alpha$ the closure of domains in$$\underline Q^{(\alpha);i}=\{Q^{(\alpha);i}_1,Q^{(\alpha);i}_2,\dots,Q^{
(\alpha);i}_j,\dots\}$$

are pairwise disjoint.\vspace{.2cm}

{\bf Property II.} The strong interiors of domains in $\underline Q^{(\alpha)}$ cover $M_\alpha$.\vspace{.2cm}

{\bf Property III.} $$Q^{(\alpha);i}_j\cap Q^{(\beta);k}_l\not=\emptyset\ \Rightarrow\ \alpha\ {\rm and}\ \beta\ {\rm can\ be\ compared\ w.r.t.
}\ \prec$$

{\noindent\it Remark:} Given a family of open subsets about strata, it is standard to impose a condition similar to Property III. (see e.g.
Mather \cite{Ma}.) It is also shown there, that arbitrary system of tubes about strata can be trimmed down so that Property III. holds.
\vspace{.3cm}

{\noindent\bf Definition 6.} A cover of each stratum by adapted charts with the above three properties is called a {\it stratified cover of}
$M$.\vspace{.4cm}

{\noindent\bf Proposition 1.} A compact $G$-manifold has a stratified cover.\vspace{.3cm}

{\noindent\it Proof:} We have already constructed a finite cover of each level-$0$ strata. For an orbit type $M_\alpha$ at the $(k+1)^{\rm
th}$ level the complement $M_\alpha\setminus\mathop{\cup}\limits_{j=1}^{j_k}U_j$ is compact so it can be covered by the strong interiors of
finitely many adapted charts, thus it is enough to choose a stratified cover for each set$$M_\alpha\cap U_j\hspace{.7cm}1\leq j\leq j_k
$$separately. For fixed $j$ domain $U_j$ is adapted to $M_\beta$ for some $\beta\prec\alpha$. Choose $\epsilon'<\epsilon$ and let $S_{\epsilon
'}\subset{\bf N}_\beta$ denote the $\epsilon'$-sphere bundle of the normal bundle of stratum $M_\beta$. The level of orbit type (exp$|_{U_j})
^{-1}(M_\alpha)\cap S_{\epsilon'}$ of the normal action is at most $k$, thus by induction we can choose a stratified cover $\underline Q'^{
(\alpha)}$ of the subset$$(S_{\epsilon'})_\alpha:=({\rm exp}|_{U_j})^{-1}(M_\alpha)\cap S_{\epsilon'}$$This means that family $\underline
Q'^{(\alpha)}$ is a finite union$$\underline Q'^{(\alpha)}=\mathop{\cup}\limits_{i=1}^k\underline Q'^{(\alpha);i}$$where each subset $
\underline Q'^{(\alpha);i}$ consists of relatively compact open sets (in the topology of $S_{\epsilon'}$) with pairwise disjoint closure.
Consider the subsets\begin{eqnarray*}{\bf R}_1&:=&\{({1\over2n+1},{2\over4n-1})\ |\ n\in{\bf Z}^+\}\\ {\bf R}_2&:=&\{({1\over2n},{4\over8n-5})
\ |\ n\in{\bf Z}^+\}\end{eqnarray*}Let\begin{eqnarray*}\underline Q_1^{(\alpha);i}&:=&\{(a,b)\times Q\ |\ (a,b)\in{\bf R}_1,\ Q\in\underline
Q'^{(\alpha);i}\}\\ \underline Q_2^{(\alpha);i}&:=&\{(a,b)\times Q\ |\ (a,b)\in{\bf R}_2,\ Q\in\underline Q'^{(\alpha);i}\}.\end{eqnarray*}
Then the exp-images of sets $\underline Q_1^{(\alpha);i},\ \underline Q_2^{(\alpha);i}\hspace{.3cm}(i=1\dots,k)$ provide a stratified cover for $M_\alpha\cap U_j$.
(Pairwise disjointness and relative compactness follow trivially; intervals $(a,b)$ can serve as new coordinates.) \hfill\rule{3mm}{3mm}\vspace{.3cm}

{\noindent\it Notation:} In each step $X$ will denote the invariant gradient-like vector field that has already been adjusted along certain
strata (so we will not re-denote $X$ in every single step).\vspace{.5cm}

For the proof of our theorem suppose inductively that relative transversality of ascending and descending submanifolds has been attained below
critical level $f(O)=c$. Fix critical orbit $O$ and a Morse-chart $\eta_O:{\bf \perp}_O(\epsilon^*)\to U_O$ about $O$. We will perturb $X$ by
an invariant vector field with support contained in the Morse-tube $U_O$ (in the case when we have more than one critical orbits at level $(f
=c)$, we choose disjoint Morse-tubes about them). This way we will not influence relative transverse intersections that have already been
established in previous steps.\vspace{.5cm}

{\noindent\it Notations:} At a point $x\in O$ let $S_x,\ S_x^-,\ S_x^+$ denote the spheres with radius $\epsilon<\epsilon^*$ (about the origin)
in subspaces $\perp_x,\ \perp^-_x,\ \perp^+_x$ respectively. Let $B=\perp_x(\epsilon^*)$ denote the $\epsilon^*$-disk about the origin of
$\perp_x$. As for the rest of this paper we will work in ball $B$, we can (ab)use notation by denoting the trajectory of vector field$$X|_B:
={\rm d}\eta_O^{-1}(X|_{U_O})={\rm grad}(f\circ\eta_O)$$through point $p\in B$ by $\lambda_p$.\vspace{.4cm}

Fix $x\in O$ and drop it from the subscripts. Let $\hat G:=G_x$ be the isotropy subgroup at $x$ and choose an invariant open neighborhood
$$\nu:=\nu_x\subset\eta_O^{-1}(f^{-1}(c-\epsilon)\cap U_O)\cap\perp_x$$of outbound shpere $S^-:=S_x^-$. Then $\nu$ is partitioned by the orbit types of
action $\hat G\times\nu\to\nu$ as$$\nu=\mathop{\cup}\limits_{[x]\preceq\alpha}\nu_\alpha$$For a subset $Z\subset\nu$ and an orbit type $[x]
\preceq\alpha$ let $Z_\alpha:=Z\cap\nu_\alpha$ be the $\alpha$-part of set $Z$ in orbit type $\nu_\alpha$. Let $\{O_1,...,O_k\}$ denote the
set of critical orbits that are connected with $O$ by a trajectory of $X$ and reside below critical level $c$. Set$$\Sigma_j:=\nu\cap\eta_O
^{-1}(W^+_{O_j}\cap U_O),\hspace{1.5cm} \Sigma:=\mathop{\cup}\limits_{j=1}^k\Sigma_j$$\vspace{.2cm}

{\noindent\it Observe:} that the intersection $W_O^-\cap W^+_j$ is relative transverse (with respect to the partition $M=\mathop{\cup}\limits
_{\alpha\in{\cal A}}M_\alpha$) if and only if for all orbit types $\alpha$ preceeded by $[x]$ the intersection $S^-_\alpha\cap(\Sigma_j)
_\alpha$ is transverse in submanifold $\nu_\alpha$. Recall notation $U_x:=\eta_O(B)$ and choose a point $p\in U_x\cap W_O^-\cap W^+_j\cap
M_\alpha$. We have\begin{eqnarray}T_p(W_O^-)_\alpha+T_p(W^+_j)_\alpha=T_p(W_O^-\cap U_x)_\alpha+T_p(W^+_j\cap U_x)_\alpha+T_pG(p)\end{eqnarray}
thus $T_p(W_O^-)_\alpha+T_p(W^+_j)_\alpha=T_pM_\alpha$ if and only if$$T_p(W_O^-\cap U_x)_\alpha+T_p(W^+_j\cap U_x)_\alpha=T_p(U_x)_\alpha$$
This means that it is enough to ensure relative transversality of intersection $S^-\cap\Sigma_j$ with respect to $\nu$.\vspace{.7cm}

{\noindent\underline{Strategy of proof:}} first we will construct the perturbing vector field on ball $B$, then we extend it along the
$G$-action to a vector field on $\perp_O(\epsilon^*)$, finally we push it forward along Morse coordinate system $\eta_O:\perp_O(\epsilon^*)
\to U_O$. We will proceed by induction on the level of orbit types of action $\hat G\times S^-\to S^-$. In the sequel
"level" will always be meant in this sense. We will use the same method to perturb vector field $X$ for both the base step and the induction
step, so by induction suppose that we need to define the perturbing vector field for orbit type $\alpha$ and we are done with all orbit types
at lower levels. We can presume $G_\alpha\subset\hat G$. Then the fixed point set of orthogonal action $G_\alpha\times\perp_x\to\perp_x$ is a
linear subspace$$V=V^-\oplus V^+$$where $V^-\subset\perp_x^-,\ V^+\subset\perp_x^+$ and dim($V^-)\not=\emptyset$.\vspace{.3cm}

\noindent Note, that the normalizer $N_\alpha$ of subgroup $G_\alpha\subset\hat G$ acts on linear subspace $V$ as well as on \vspace{.2cm}

$V_\alpha^-:=[V^-$ less the points that are fixed by a group strictly larger than $G_\alpha$]\vspace{.3cm}

\noindent Let $Q^*$ be an $\alpha$-slice for action $N_\alpha\times V_\alpha^-\to V_\alpha^-$ such that it is a union of (open) rays in $V_\alpha^-$
(thus it is a cone ${\bf C}(Q^*\cap S^-)$ over $\alpha$-slice $Q^*\cap S^-$ for action $N_\alpha\times S^-\to S^-$. Let $D^+\subset V^+$ be a
small disc about the origin. Then $Q^*\oplus D^+$ is an $\alpha$-slice of action $N_\alpha\times V\to V$ and also for action $\hat G\times
\perp_x\to\perp_x$, thus$$(Q^*\oplus D^+)\cap\nu$$is an $\alpha$-slice for action $\hat G\times\nu\to\nu$. This leads to the following
conclusion, which is crucial for the proof:\vspace{.3cm}

\noindent{\it Observation 1.:} For an $\alpha$-slice $Q^{**}$ of action $\hat G\times S^-\to S^-$ there exists an $\alpha$-slice of action$$
\hat G\times\nu\to\nu$$which is a product of $Q^{**}$ and the fiber. We will call such a slice a {\it product slice associated} to $Q^{**}$
and its image under $G$ (i.e. the union of paths through its points) is the {\it G-extension} of the associated slice.\vspace{.3cm}

Level-$0$ orbit types are compact, so are their projections to orbit space $\cal M$. In the base step of induction we define the perturbing
vector fields along these projections and then we will lift them into invariant vector fields defined on $M$. This is done in complete analogy with the non-equivariant case (see Smale \cite{S}).\vspace{.3cm}

\noindent{\it Remark:} We could have chosen to work out the whole reasoning in the orbit space, lifting the resulting vector fields
(isotopies) to $M$ afterwards. In spite of the simplicity of some parts of this path of reasoning, other technical difficulties would have
arisen (e.g. the quotient space $\cal V$ is not a linear space, etc).\vspace{.4cm}

{\noindent\bf Construction of the stratified cover:} first we choose a family of $\alpha$-slices$$Q^{(\alpha);i*}_j\subset S^-_\alpha\cap V^-
\hspace{.8cm}(i=1,\dots,k_\alpha,\ j=1,\dots,n,\dots)$$for the action $$\hat G\times S^-\to S^-$$together with associated product slices
$$Q^{(\alpha);i*}_j\oplus D^+$$so that $$\overline{\hat GQ^{(\alpha);i*}_j}\cap\overline{\hat GQ^{(\alpha);i*}_{j'}}=\emptyset\hspace{.4cm}
(j\not=j')$$Intersections$$\nu\cap{\bf C}(Q^{(\alpha);i*}_j\oplus D^+)$$will be $\alpha$-slices for action$$\hat G\times\nu\to\nu$$so
 if ${\bf n}_\alpha$ is the normal bundle of $\nu_\alpha\subset\nu$ then with the aid of a small disc bundle of$${\bf n}_\alpha|_{\hat G{\bf
C}(Q^{(\alpha);i*}_j\oplus D^+)\cap\nu}$$we can define a stratified cover$$\underline Q^{(\alpha)}=\underline Q^{(\alpha);1}\cup...\cup
\underline Q^{(\alpha);k_\alpha}$$of the $\alpha$-part $\nu_\alpha$ (the cover is taken in $\nu$). Choose a smaller disk $\overline{D'^+}
\subset D^+$ and fix compact subsets$$K^{(\alpha);i*}_j\subset Q^{(\alpha);i*}_j,\hspace{1cm}K^{(\alpha);i}_j\subset Q^{(\alpha);i}_j$$With
$\epsilon'<\epsilon$ define\begin{eqnarray}\tilde Q^{(\alpha);i}_j:=\{q\in\lambda_p\ |\ p\in Q^{(\alpha);i}_j\ {\rm and}\ c-\epsilon
<f(q)<c-\epsilon'\}\\ \tilde K^{(\alpha);i}_j:=\{q\in\lambda_p\ |\ p\in K^{(\alpha);i}_j\ {\rm and}\ c-\epsilon<f(q)<c-\epsilon'\}\end
{eqnarray}

The $\alpha^{\rm th}$ step will consist of $k_\alpha$ substeps so that in the $i^{\rm th}$ substep the support of the perturbing vector field
is contained in the extension$$G\left(\mathop{\cup}\limits_{j=1}^\infty\tilde Q^{(\alpha);i}_j\right)\hspace{.3cm}(i=1,\dots,k
_\alpha)$$\vspace{.3cm}

\noindent The following observations can be made:\vspace{.4cm}

{\noindent\it Observation 2.} The perturbed vector field remains the same on lower level strata, so does  the transversality of intersections
$S^-_\beta\cap\Sigma_\beta\hspace{.2cm}(\beta\prec\alpha)$ (that have already been established by induction).\vspace{.2cm}

{\noindent\it Observation 3.} By the fact that transverse intersections are stable under small perturbations of class $C^1$, in each sub-step
we can choose the perturbing vector field to be so small that it will not destroy transverse intersections that had already been established
in previous sub-steps. Thus it is enough to describe the $i^{\rm th}$ sub-step, or, as the domains$$\{\tilde Q^{(\alpha);i}_j\ |\ j=1,\dots,n,\dots\}$$are
pairwise disjoint, it is enough to describe how to modify $X$ within one such domain.\vspace{.2cm}

{\noindent\it Observation 4.} The bound we imposed in the Theorem divided by the total number of sub-steps (i.e.$$\varepsilon':={\varepsilon\over
\mathop{\sum}\limits_{\alpha\in{\cal A}}k_\alpha})$$provides a bound we should use in formula (4) in each sub-step.\vspace{.2cm}

{\noindent\it Observation 5.} For each triple $(\alpha,i,j)$ there is a bound $\varepsilon^{(\alpha);i}_j$ so that if in the $i^{\rm th}$
sub-step we choose the perturbing vector field with $C^1$-norm measured on domain $\tilde Q^{(\alpha);i}_j$ smaller than $\varepsilon^{(\alpha
);i}_j$, then its $C^1$-norm defined by formula (4) is smaller than $\varepsilon'$ (this is a well-known fact that follows from relative
compactness of the domains, see e.g. Sternberg \cite{St}.)\vspace{.2cm}

{\noindent\it Observation 6.} The major difficulty is to ensure that after the final step we end up with a $C^1$ vector field. This
will give us additional conditions on the size of perturbation we can make in each sub-step. These conditions are described in the Lemma
below.\vspace{.5cm}

\noindent By the above observations it is enough to show that for arbitrary indices $(\alpha,i,j)$ (that we fix and drop, using notation
$\varepsilon":=\varepsilon^{(\alpha);i}_j$) the following holds:
\vspace{.3cm}

{\noindent\bf Proposition 2.} Given a domain $\tilde Q=\tilde Q^{(\alpha);i}_j$ with compact subset $\tilde K=\tilde K^{(\alpha);i}_j$ and $\varepsilon">0$
there exists an invariant vector field $Y$ with support in $G\tilde Q$ such that:\vspace{.3cm}

(i) Vector field $X+$ d$\eta_O(Y)$ is also gradient like for $G$-Morse function $f$.\vspace{.2cm}

(ii) $\|Y\|_{1;\tilde Q}<\varepsilon"$ (i.e. the $C^1$-norm on $G\tilde Q$ is smaller than $\varepsilon"$.)\vspace{.2cm}

(iii) The intersection$$[S^-_\alpha\cap\hat GK^{(\alpha);i}_j]\cap[\Sigma_\alpha\cap\hat GK^{(\alpha);i}_j]$$

(notations stand for objects of the flow of vector field $X|_B+Y$) is transverse

in orbit type $\nu_\alpha$, thus by formula (5) intersection\begin{eqnarray}[W_O^-\cap G\tilde K\cap M_\alpha]\cap[W_{O'}^+\cap G\tilde K\cap M_\alpha]
\end{eqnarray}

is transverse in orbit type $M_\alpha$ for any other critical orbit $O'$.\vspace{.3cm}

{\noindent\it Proof:} dim$(V^+)=0$ implies $\nu_\alpha=S^-_\alpha$ thus transversality of intersections in formula (8) follow trivially.

Otherwise note that $V\cap\Sigma_j$ is the fixed point set of subgroup $G_\alpha\subset\hat G$ for the action $\hat G\times\Sigma_j\to\Sigma_j
$, thus it is a submanifold of $\Sigma_j$. Let $Q_p^*$ be an $\alpha$-slice at point $p\in\Sigma_j\cap V_\alpha^-(=\Sigma_j\cap S^-_\alpha)$
for action $N_\alpha\times(V\cap\Sigma_j)\to(V\cap\Sigma_j)$. Then $V_p^*:=T_pQ_p^*$ is an affine subspace of $V$, so it splits into the sum$$V_p^*=V_p^{*-}\oplus V_p
^{*+}$$where $V_p^{*-}:=V_p^*\cap V^-$. Intersection $(\Sigma_j)_\alpha\cap S^-_\alpha$ is $\alpha$-transverse iff dim$(V_p^{*+})$=dim$(V^+)$,
or in other words if the origin is a regular value of the projection to the second factor $$P^+|_{V\cap\Sigma_j}:V\cap\Sigma_j\to V^+$$(see
also formula (5)).

By Sard's theorem the set of critical values of the above projection is of measure-$0$ for $j=1,\dots,k$ thus the same holds for their union.
This means that we can choose an arbitrarily small vector ${\bf v}\in D^+$ so that the constant section ${\bf C}(K^{(\alpha);i*}_j\oplus {\bf
 v})\cap\nu_\alpha$ intersects $(\Sigma_j)_\alpha$ $\alpha$-transversely (i.e., relative to $\nu_\alpha$) for $j=1,\dots,k$.\vspace{.3cm}

\noindent Choose cutoff function$$\phi:Q^{(\alpha);i*}_j\to[0,1],\hspace{.8cm}\phi(K^{(\alpha);i*}_j)=1,\hspace{.8cm}{\rm supp}(\phi)\subset
 Q^{(\alpha);i*}_j$$and let $H_D$ be an isotopy of disk $D^+$ which moves the origin into $\bf v$ and fixes a neighborhood of the boundary.
Define isotopy$$H:(Q^{(\alpha);i*}_j\oplus D^+)\times[0,1]\to Q^{(\alpha);i*}_j\oplus D^+,\hspace{.8cm}H((p,{\bf w}),t)=(p,H_D({\bf w},\phi
(p)t)$$Using rays this defines an isotopy of set ${\bf C}(Q^{(\alpha);i*}_j\oplus D^+)\cap\nu$. An application of the argument in Milnor
(\cite{Mi}, pp. 42-43) to the restriction $X_B|_F$ where$$F:=\{q\in\lambda_p\ |\ p\in{\bf C}(Q^{(\alpha);i*}_j\oplus D^+)\cap\nu\ {\rm and}\
c-\epsilon<f(q)<c-\epsilon'\}$$produces a vector field $Y^*$ which is tangent to $F$ so that the difference between moving along the flow-lines of vector fields
 $X_B|_F$ and $X_B|_F+Y^*$ shows up in the application of map $H_1$ on level set $(f=c-\epsilon)$. Let $\hat Y$ be the extension of $Y^*$ along the
action onto set $\hat GF$. This set is an open subset of the $\alpha$-part $B_\alpha$ of disk $B=\perp_x(\epsilon^*)$.

\noindent For a $\delta>0$ choose cutoff function$$\psi:[0,\delta]\to[0,1],\hspace{.8cm}\psi(0)=1,\ \psi(\delta)=(0),\ {{\rm d}\psi\over{\rm d}t}\leq0.
$$Let ${\bf N}_\alpha$ be the normal bundle of the $\alpha$-part $B_\alpha\subset B$ (taken in ball $B$). Consider the Levi-Civita connection
of the Euclidean metric on vector bundle ${\bf N}_\alpha$ and let $Y'$ be the horizontal lift of $\hat Y$ (i.e. a vector field on ${\bf N}
_\alpha$ with vertical components $=0$ at each point so that d$\pi(Y'|_p)=\hat Y|_{\pi(p)}$). Let$$Y={\rm exp}_*(\psi(\|v\|)Y'|_v)\hspace
{.8cm}(v\in{\bf N}_\alpha)$$ and extend $Y$ via the action onto tube $U_O$. Let $Y$ be the $0$ vector field outside of set $G{\rm exp}({\bf
N}_\alpha(\delta)|_{\hat GF})$. Then clauses (i) and (iii) hold for vector field $X_B+Y$ (its restriction to a smaller neighborhood of critical
orbit $O$ is gradient-like for $f$).

As in Definition 5., choose a coordinate system $(x^2,\dots,x^{m_\alpha})$ on $\alpha$-slice

\noindent $Q^{(\alpha);i*}_j\oplus D^+$. Use $G$-Morse function $f$ as the first coordinate $x^1$. As in clause (iv) of Definition 5., fixing an
orthonormal frame bundle on ${\bf N}_\alpha|_{F}$ introduces further coordinates $x^{m_\alpha+1},\dots,x^{m_\alpha+m_\alpha^\perp}$.
Finally, supplementing these with the coordinates in formula (2) provides a finite family of coordinate charts on invariant open set $G\tilde
Q$. Note, that only the first $m_\alpha+m_\alpha^\perp$ coordinates of vector field $Y$ are non-zero and the coordinates themselves are bounded by $\|{\bf v}
\|$. By relative compactness of set $Q^{(\alpha);i*}_j$ we can choose an upper bound $N$ for the absolute value of the derivatives of cutoff functions
$\phi$ and $\psi$. Then $$2\left(N\|{\bf v}\|+{2\|{\bf v}\|\over\epsilon-\epsilon'}\right)$$will serve as an upper bound for the absolute
values of the first derivatives of coordinates of vector field $Y$. This shows that by choosing $\|{\bf v}\|$ small enough, one can arrange
that clause (ii) of our Proposition holds as well.\hfill\rule{3mm}{3mm}\vspace{.3cm}

Differentiability of the final vector field can be ensured by choosing $\|{\bf v}\|$ according to how close the domain  $\tilde Q^{(\alpha);
i}_j$ is from the boundary of $B_\alpha$. Let$${\rm Fr}(B_\alpha):=\overline B_\alpha\setminus B_\alpha$$be the {\it frontier} of orbit type
$B_\alpha$ and let $d$ stand for the Euclidean distance on ball $B$. Then, by relative compactness of set $Q^{(\alpha);i*}_j$ in $S_\alpha$,
the distance $$d_{ij}:=d(\hat G\tilde Q^{(\alpha);i}_j,{\rm Fr}(B_\alpha))>0$$\vspace{.2cm}

\noindent In the $i^{\rm th}$ substep the components of supp$(Y)$ are contained in subset$$\mathop{\cup}\limits_{j=1}^\infty G\tilde Q^{(
\alpha);i}_j$$By invariance it is enough to prove differentiability of restriction $Y|_B$. By construction, $Y$ is a smooth vector field in a
tube about $\alpha$-part $B_\alpha$ (i.e. on a set$${\rm exp}({\bf N}_\alpha(\zeta))$$ where ${\bf N}_\alpha(\zeta)$ is an appropriate disk
bundle of the normal bundle of $B_\alpha\subset B$) and $Y$ is the zero vector field outside of this tube. This shows that we can have problems
with differentiability of vector field $Y$ only at the points of the frontier Fr$(B_\alpha)$. We will prove that imposing an additional
condition ensures that vector field $Y$ is of class $C^1$ at points of Fr$(B_\alpha)$. We will formulate this condition by utilizing the fact that
$B$ is a Euclidean disk, thus its tangent bundle is trivial by natural identification$$TB=B\times{\bf R}^n$$($n$=dim$(B)$.) This way we can
look at vector field $Y|_B$ as a map$$Y_B:B\to{\bf R}^n$$ \vspace{.2cm}

{\noindent\bf Lemma.} Vector field $Y$ is of class $C^1$ whenever$$\|Y(q)\|<d_{ij}^3\hspace{.8cm}(q\in \tilde Q^{(\alpha);i}_j,\ i=1,\dots,
k_\alpha,\ j\in{\bf N})$$

{\it Proof} A point $q\in$Fr$(B_\alpha)$ belongs to $q\in B_\beta$ for some orbit type $\beta\prec\alpha$. By definition $Y$ is differentiable
at point $q$ if\begin{eqnarray}\mathop{\lim}\limits_{p\to q}{\|Y(p)-Y(q)\|\over\|p-q\|}=0\end{eqnarray}Value $Y(q)={\bf 0}$ for points of the frontier, thus $Y(p)\not={\bf
0}$ implies that $p\in\hat G\tilde Q^{(\alpha);i}_j$ for some indices $i\in\{1,\dots,k_\alpha\},\ j\in{\bf N}$. But then$$\|Y(p)\|<d_{ij}^3
\leq\|p-q\|^3$$thus the limit in formula (9) is indeed $0$, together with the limit of the first partial derivatives of $Y$.\hfill\rule{3mm}{3mm}\vspace{.3cm}

{\noindent\it Remark:} The proof of the lemma has been built heavily on the Eucledian structure on disk $B$. In a general setting a somewhat
more complicated proof would work: one considers the isotopy induced by time-dependent vector field$$Y_t:=\xi(t)Y\hspace{.4cm}(t\in[0,1])
$$where $\xi:[0,1]\to[0,1],\ \xi[0,\epsilon)=0$ is some cutoff function. One puts the analoguos condition on the displacement of this isotopy
(measured in the metric distance on $B$). Then a similar argument proves that the isotopy is of class $C^1$, consequently vector field $Y$ is
also of class latex GSM.txt$C^1$.\vspace{.3cm}

This paper is an improved version of Chapter 3. of my Ph.D. thesis \cite{M}. I wish to express my gratitude to my advisor, Professor Dan Burghelea for
his help.

Imre Major

Gabor Denes School of Informatics and

Central European University, Budapest

\end{document}